\documentclass[global]{svjour}
\usepackage{amssymb}
\usepackage{amsfonts}
\usepackage{amsmath}

\input{tcilatex}

\begin{document}

\title{Fixed point theorems for nonconvex valued correspondences and
applications in game theory}
\author{Monica Patriche}
\institute{University of Bucharest 
%TCIMACRO{\TeXButton{email}{\email{monica.patriche@yahoo.com}}}%
%BeginExpansion
\email{monica.patriche@yahoo.com}%
%EndExpansion
}
\mail{University of Bucharest, Faculty of Mathematics and Computer Science, 14
Academiei Street, 010014 Bucharest, Romania}
\maketitle

\begin{abstract}
{\small In this paper, we introduce several types of correspondences: weakly
naturally quasiconvex, *-weakly naturally quasiconvex, weakly biconvex and
correspondences with *--weakly convex graph and we prove some fixed point
theorems for these kinds of correspondences. As a consequence, using a
version of W. K. Kim's quasi-point theorem, we obtain the existence of
equilibria for a quasi-game. }
\end{abstract}



\section{Introduction}

The aim of this paper is to prove some fixed points theorems for
correspondences which are not continuous or convex valued and to give
applications in game theory.

The significance of equilibrium theory stems from the fact that it develops
important tools (as fixed point and selection theorems) to prove the
existence of equilibrium for different types of games. In 1950, J. F. Nash
[15] first proved a theorem of equilibrium existence for games where the
player's preferences were representable by continuous quasi-concave
utilities. G. Debreu's works on the existence of equilibrium in a
generalized N-person game or on an abstract economy [6] were extended by
several authors. In [16] W. Shafer and H. Sonnenschein proved the existence
of equilibrium of an economy with finite dimensional commodity space and
irreflexive preferences represented as correspondences with open graph. N.
C. Yannelis and N. D. Prahbakar [19] developed new techniques based on
selection theorems and fixed-point theorems. Their main result concerns the
existence of equilibrium when the constraint and preference correspondences
have open lower sections. They worked within different frameworks (countable
infinite number of agents, infinite dimensional strategy spaces). K. J.
Arrow and G. Debreu proved the existence of Walrasian equilibrium in [3]. In
[20], X. Z. Yuan proposed a model of abstract economy more general than that
introduced by Borglin and Keing in [4].

Within the last years, a lot of authors generalized the classical model of
abstract economy. For example, K. Vind [18] defined the social system with
coordination, X. Z. Yuan [20] proposed the model of the general abstract
economy. Motivated by the fact that any preference of a real agent could be
unstable by the fuzziness of consumers' behaviour or market situation, W.
K.\ Kim and K. K.\ Tan [12] defined the generalized abstract economies. Also
W. K. Kim [13] obtained a generalization of the quasi fixed-point theorem
due to I. Lefebvre [14], and as an application, he proved an existence
theorem of equilibrium for a generalized quasi-game with infinite number of
agents. W. K. Kim's result concerns generalized quasi-games where the
strategy sets are metrizable subsets in locally convex linear topological
spaces.

Biconvexity was studied by R. Aumann, S. Hart in [2] and J. Gorski, F.
Pfeuffer and K. Klamroth in [10].

An open problem of the fixed point theory is to prove the existence of fixed
points for correspondences without continuity or convex values. X. Ding and
He Yiran introduced in [7] the correspondences with weakly convex graph to
prove a fixed point theorem. The result concerning the existence of the
affine selection on a special type of sets (simplex) proves to be redundant,
since a correspondence $T:X\rightarrow 2^{Y}$ has an affine selection if and
only if it has a weakly convex graph. This result is stronger than it needs
in order to obtain a fixed point theorem. We try to weaken these conditions
by defining several types of correspondences which are not continuous or
convex valued: weakly naturally quasiconvex, *-weakly naturally quasiconvex,
correspondences with *--weakly convex graph and weakly biconvex
correspondences. We prove fixed point theorems for these kinds of
correspondences and using a version of W. K. Kim's quasi-point theorem, we
prove the existence of equilibria for a quasi-game. We use the continuous
selection technique introduced by N. C. Yannelis and N. D. Prahbakar in [19].

The paper is organized in the following way: Section 2 contains
preliminaries and notation. The fixed point theorems are given in Section 3
and the equilibrium theorems are stated in Section 4.\bigskip

\section{\textbf{Preliminaries and notation\protect\smallskip \protect%
\medskip }}

Throughout this paper, we shall use the following notations and definitions:

Let $A$ be a subset of a topological space $X$.

1. 2$^{A}$ denotes the family of all subsets of $A$.

2. cl $A$ denotes the closure of $A$ in $X$.

3. If $A$ is a subset of a vector space, co$A$ denotes the convex hull of $A$%
.

4. If $F$, $T:$ $A\rightarrow 2^{X}$ are correspondences, then co$T$, cl $T$%
, $T\cap F$ $:$ $A\rightarrow 2^{X}$ are correspondences defined by $($co$%
T)(x)=$co$T(x)$, $($cl$T)(x)=$cl$T(x)$ and $(T\cap F)(x)=T(x)\cap F(x)$ for
each $x\in A$, respectively.

5. The graph of $T:X\rightarrow 2^{Y}$ is the set Gr$(T)=\{(x,y)\in X\times
Y\mid y\in T(x)\}$

6. The correspondence $\overline{T}$ is defined by $\overline{T}(x)=\{y\in
Y:(x,y)\in $cl$_{X\times Y}$Gr$T\}$ (the set cl$_{X\times Y}$Gr$(T)$ is
called the adherence of the graph of T)$.$

It is easy to see that cl$T(x)\subset \overline{T}(x)$ for each $x\in
X.\medskip $

\begin{definition}
Let $X$, $Y$ be topological spaces and $T:X\rightarrow 2^{Y}$ be a
correspondence. $T$ is said to be \textit{upper semicontinuous} if for each $%
x\in X$ and each open set $V$ in $Y$ with $T(x)\subset V$, there exists an
open neighborhood $U$ of $x$ in $X$ such that $T(x)\subset V$ for each $y\in
U$.\medskip
\end{definition}

Let $X\subset E_{1}$ and $Y\subset E_{2}$ be two nonempty convex sets, $%
E_{1},E_{2}$ be topological vector spaces and let $B\subset X\times Y.$

\begin{definition}[2]
The set $B\subset X\times Y$ is called a \textit{biconvex set} on $X\times Y$
if the section $B_{x}=\left\{ y\in Y:(x,y)\in B\right\} $ is convex for
every $x\in X$ and the section $B_{y}=\left\{ x\in X:(x,y)\in B\right\} $ is
convex for every $y\in Y.\medskip $
\end{definition}

\begin{definition}[2]
Let $(x_{i},y_{i})\in X\times Y$ for $i=1,2,...n.$ A convex combination $%
(x,y)=\tsum\limits_{i=1}^{n}\lambda _{i}(x_{i},y_{i})$, (with $%
\tsum\limits_{i=1}^{n}\lambda _{i}=1,$ $\lambda _{i}\geq 0$ $i=1,2,...,n$)
is called \textit{biconvex combination} if $x_{1}=x_{2}=...=x_{n}=x$ or $%
y_{1}=y_{2}=...=y_{n}=y.\medskip $
\end{definition}

\begin{theorem}
(Aumann and Hart [2]). A set $B\subseteq X\times Y$ is biconvex if and only
if $B$ contains all biconvex combinations of its elements.\medskip
\end{theorem}

\begin{definition}[2]
Let $A\subseteq X\times Y$ be a given set. The set $H:=\tbigcap
\{A_{I}:A\subseteq A_{I},$
\end{definition}

\noindent $A_{I}$ is biconvex\} is called \textit{biconvex hull of }$A$ and
is denoted biconv$(A).\medskip $

\begin{theorem}
(Aumann and Hart [2]). The \textit{biconvex hull of a set} $A$ is biconvex.
Furthermore, it is the smallest biconvex set (in the sens of set inclusion),
which contains $A.\medskip $
\end{theorem}

\begin{lemma}
(Gorski, Pfeuffer and Klamroth [10]). Let $A\subseteq X\times Y$ be a given
set. Then biconv$(A)\subseteq $conv$(A).\medskip $
\end{lemma}

\textit{NOTATION.} We denote the standard $(n-1)$- dimensional simplex by $%
\Delta _{n-1}=\{(\lambda _{1},\lambda _{2},...,\lambda _{n})\in \mathbb{R}%
^{n}:\overset{n}{\underset{i=1}{\tsum }}\lambda _{i}=1,\lambda _{i}\geqslant
0,i=1,2,...,n\}$.

\section{Selection theorems and fixed point theorems}

An open problem of the fixed point theory is to prove the existence of fixed
points for correspondences without continuity or convex values. In this
section we introduce some types of correspondences which are not continuous
or convex valued and prove selection theorems and fixed point theorems.

First, we introduce the concept of weakly naturally quasiconvex
correspondence.

\begin{definition}
Let $X,Y$ be nonempty convex subsets of topological vector spaces $E,$
respectively\ $F$.\textit{\ }The correspondence\textit{\ }$%
T:X\longrightarrow 2^{Y}$ is said to be \emph{weakly naturally quasiconvex
(WNQ) } if for each $n\in \mathbb{N}$ and for each finite set $%
\{x_{1},x_{2},...,x_{n}\}\subset X$, there exists $y_{i}\in T(x_{i})$ , $%
(i=1,2,...,n)$ and $g=(g_{1},g_{2},...,g_{n}):\Delta _{n-1}\rightarrow
\Delta _{n-1}$ a bijective function depending on $x_{1},x_{2},...,x_{n}$
with $g_{i}$ continuous, $g_{i}(1)=1,$ $g_{i}(0)=0$ for each\textit{\ }$%
i=1,2,...n$, such that for every $(\lambda _{1},\lambda _{2},...,\lambda
_{n})\in \Delta _{n-1}$, there exists $y\in T(\overset{n}{\underset{i=1}{%
\tsum }}\lambda _{i}x_{i})$ and $y=\overset{n}{\underset{i=1}{\tsum }}%
g_{i}(\lambda _{i})y_{i}.\medskip $
\end{definition}

\begin{remark}
A weakly naturally quasiconvex correspondence may not be continuous or
convex valued.
\end{remark}

We give an economic interpretation of the weakly naturally quasiconvex
correspondences.

We consider an abstract economy $\Gamma =(X_{i},A_{i},P_{i})_{i\in I}$ with $%
I$ - the set of agents. Each agent can choose a strategy from the set $X_{i}$
and has a preferrence correspondence $P_{i}:X=\tprod\limits_{i\in
I}X_{i}\rightarrow 2^{X_{i}}$ and a constraint correspondence $%
A_{i}:X=\tprod\limits_{i\in I}X_{i}\rightarrow 2^{X_{i}}.$ The traditional
approach considers that the preferrence of agent $i$ is characterized by a
binary relation $\succeq _{i}$ on the set $X_{i}.$ A real valued function $%
u_{i}:X\rightarrow \mathbb{R}$ that satisfies $x\succeq _{i}y$ $%
\Leftrightarrow $ $u_{i}(x)\geq u_{i}(y)$ is called an utility function of
the preferrence $\succeq _{i}.$ The relation between the utility function $%
u_{i}$ and the preferrence correspondence $P_{i}$, for each agent $i$ is:

$P_{i}(x)=\left\{ y_{i}\in X_{i}:u_{i}(x,y_{i})>u_{i}(x,x_{i})\right\} ,$
where, in this case, $u_{i}:X\times X_{i}\rightarrow \mathbb{R}.$

The aim of the equilibrium theory is to maximize each agent's utility on a
strategy set.

For the case that, for each index $i$, $P_{i}$ is a weakly naturally
quasiconvex correspondence, the interpretation is the following: for all
certain amounts $x^{1},x^{2},...x^{n}\in X,$ the agent $i$ with the
correspondence $P_{i}$ will always prefer $y_{i},$ the weighted average of
some quantities $y_{i}^{k}\in P_{i}(x^{k}),$ $k=1,n.$ This implies that
there exist $y_{i}^{1}\in P_{i}(x^{1}),y_{i}^{2}\in
P_{i}(x^{2}),...,y_{i}^{n}\in P_{i}(x^{n})$ and $g=(g_{1},g_{2},...,g_{n}):%
\Delta _{n-1}\rightarrow \Delta _{n-1}$ a bijective function depending on $%
x^{1},x^{2},...,x^{n}$ with $g_{i}$ continuous for each\textit{\ }$%
i=1,2,...,n,$ such that, for each $\lambda \in \Delta _{n-1},$ there exists $%
y_{i}=\overset{n}{\underset{k=1}{\tsum }}g_{i}(\lambda _{k})y_{i}^{k}$ and $%
y_{i}\in P_{i}(\overset{n}{\underset{k=1}{\tsum }}\lambda _{k}x^{k})$ (i.e.,
if there exist utility functions $u_{i}:X\times X_{i}\rightarrow \mathbb{R}$
such that, if $u_{i}(x^{k},y_{i}^{k})>u_{i}(x^{k},x_{i}^{k})$ for every $%
k\in \{1,2,...n\},$ we have that $y_{i}\in A_{i}(\overset{n}{\underset{k=1}{%
\tsum }}\lambda _{k}x^{k})$ and $u_{i}(\overset{n}{\underset{k=1}{\tsum }}%
\lambda _{k}x^{k},y_{i})>u_{i}(\overset{n}{\underset{k=1}{\tsum }}\lambda
_{k}x^{k},(\overset{n}{\underset{k=1}{\tsum }}\lambda _{k}x^{k})_{i}).$

\begin{theorem}
(selection theorem). \textit{Let }$K$\textit{\ be a simplex in a topological
vector space }$F$ and $Y$\textit{\ be a non-empty convex subset of a
topological vector space }$E.$ \textit{Let }$T:K\rightarrow 2^{Y}$\textit{\
be a weakly naturally quasiconvex correspondence. Then, }$T$\textit{\ has a
continuous selection on }$K$\textit{.\medskip }
\end{theorem}

\textit{Proof.} Assume that $K$ is a simplex, i.e., the convex hull of an
affinely independent set $\{a_{1},a_{2},...,a_{n}\}.$ Since $T$ is weakly
naturally quasiconvex, there exist $b_{i}\in T(a_{i})$, $(i=1,2,...,n)$ and $%
g=(g_{1},g_{2},...,g_{n}):\Delta _{n-1}\rightarrow \Delta _{n-1}$ a
bijective function with $g_{i}$ continuous for each\textit{\ }$i=1,2,...,n$,
such that for every $(\lambda _{1},\lambda _{2},...,\lambda _{n})\in \Delta
_{n-1}$, there exists $y\in T(\overset{n}{\underset{i=1}{\tsum }}\lambda
_{i}a_{i})$ with $y=\overset{n}{\underset{i=1}{\tsum }}g_{i}(\lambda
_{i})y_{i}.$

Since $K$ is a $(n-1)$-dimensional simplex with the vertices $%
a_{1},...,a_{n},$ there exists unique continuous functions $\lambda
_{i}:K\rightarrow \mathbb{R},$ $i=1,2,...,n$ such that for each $x\in K,$ we
have $(\lambda _{1}(x),\lambda _{2}(x),...,\lambda _{n}(x))\in \Delta _{n-1}$
and $x=\overset{n}{\underset{i=1}{\tsum }}\lambda _{i}(x)a_{i}.$

Let's define $f:K\rightarrow Y$ by

$f(a_{i})=b_{i}$ $(i=1,...,n)$ and

$f(\overset{n}{\underset{i=1}{\tsum }}\lambda _{i}a_{i})=\overset{n}{%
\underset{i=1}{\tsum }}g_{i}(\lambda _{i})b_{i}\in T(x).$

We show that $f$ is continuous.

Let $(x_{m})_{m\in N}$ be a sequence which converges to $x_{0}\in K,$ where $%
x_{m}=\overset{n}{\underset{i=1}{\tsum }}\lambda _{i}(x_{m})a_{i}$ and $%
x_{0}=$ $\overset{n}{\underset{i=1}{\tsum }}\lambda _{i}(x_{0})a_{i}.$ By
the continuity of $\lambda _{i},$ it follows that for each $i=1,2,...,n$, $%
\lambda _{i}(x_{m})\rightarrow \lambda _{i}(x_{0})$ as $m\rightarrow \infty
. $ Since $\ g_{1},g_{2},...,g_{n}$ are continuous, we have $g_{i}(\lambda
_{i}(x_{m}))\rightarrow g_{i}(\lambda _{i}(x_{0}))$ as $m\rightarrow \infty
. $ Hence, $f(x_{m})\rightarrow f(x_{0})$ as $m\rightarrow \infty ,$ i.e. $f$
is continuous.

We proved that $T$ has a continuous selection on $K.\medskip $

By Brouwer's fixed point theorem, we obtain the following fixed point
theorem for weakly naturally quasiconvex correspondences.

\begin{theorem}
\textit{Let }$K$\textit{\ be a simplex in a topological vector space }$F.$ 
\textit{Let }$T:K\rightarrow 2^{K}$\textit{\ be a weakly naturally
quasiconvex correspondence. Then, }$T$\textit{\ has a fixed point in }$K$%
\textit{.}
\end{theorem}

\textit{Proof.} By Theorem 3, $T$ has a continuous selection on $K,$ $%
f:K\rightarrow K.$

Since $f$ has a fixed point $x^{\ast }\in K,$ we have that $x^{\ast
}=f(x^{\ast })\in T(x^{\ast }).$\textit{\medskip }

\textit{NOTATION.} For the correspondence $T:X\rightarrow 2^{Y}$ and for the
set $V\in Y,$ we denote $T_{V}$ the correspondence $T_{V}:X\rightarrow 2^{Y}$%
, defined by $T_{V}(x)=(T(x)+V)\cap Y$ for each $x\in X.$

If $Y=K,$ we obtain the following fixed point theorem:

\begin{theorem}
Let $K$ be a \textit{simplex in} a topological vector space $F$ and let $%
T:K\rightarrow 2^{K}$ be a correspondence. Assume that for each neighborhood 
$V$ of the origin in $F$, there is $T^{V}:K\rightarrow 2^{K}$ a weakly
naturally quasiconvex correspondence such that Gr$T^{V}\subset $clGr$T_{V}$.
Then there exists a point $x^{\ast }\in K$ such that $x^{\ast }\in \overline{%
T}(x^{\ast }).$
\end{theorem}

To prove Theorem 5, we need the following lemma from [20].

\begin{lemma}[20]
Let $X$ be a topological space, $Y$ be a non-empty subset of a topological
vector space E, \ss\ be a base of the neighborhoods of $0$ in $E$ and $%
T:X\rightarrow 2^{Y}.$ If $x^{\ast }\in X$ and $\widehat{y}\in Y$ are such
that $\widehat{y}\in \cap _{V\in \text{\ss }}\overline{T_{V}}(x^{\ast }),$
then $\widehat{y}\in \overline{T}(x^{\ast }),$ where $\overline{T}$ $%
:X\rightarrow 2^{Y}$ is defined by $\overline{T}(x)=\{y\in Y:(x,y)\in $cl$%
_{X\times Y}$Gr$T\}.\medskip $
\end{lemma}

\textit{Proof of Theorem 5.} Let \ss\ denote the family of all neighborhoods
of zero in $F.$ Let $V\in $\ss $.$ By the fixed point theorem 4, it follows
that for each neighborhood $V$\ of the origin in $Y,$ there exists $%
x_{V}^{\ast }\in T^{V}(x_{V}^{\ast })\subset (T(x_{V}^{\ast })+V)\cap K.$

For each $V\in \text{\ss ,\ we \ define }Q_{V}=\{x\in K:$ $x\in (T(x)+V)\cap
K\}.$

$Q_{V}$ is nonempty since $x_{V}^{\ast }\in Q_{V},$ then cl$Q_{V}$ is
nonempty.

We prove that the family $\{$cl$Q_{V}:V\in \text{\ss }\}$ has the finite
intersection property.

Let $\{V^{(1)},V^{(2)},...V^{(n)}\}$ be any finite set of $\text{\ss }.$ Let 
$V=\underset{k=1}{\overset{n}{\cap }}V^{(k)}$, then $V\in \text{\ss }$.
Clearly $Q_{V}\subset \underset{k=1}{\overset{n}{\cap }}Q_{V^{(k)}}$ so that 
$\underset{k=1}{\overset{n}{\cap }}Q_{V^{(k)}}\neq \emptyset .$ Then $%
\underset{k=1}{\overset{n}{\cap }}$cl$Q_{V^{(k)}}\neq \emptyset .$

Since $K$ is compact and the family $\{$cl$Q_{V}:V\in \text{\ss }\}$ has the
finite intersection property, we have that $\cap \{$cl$Q_{V}:V\in \text{\ss }%
\}\neq \emptyset .$ Take any $x^{\ast }\in \cap \{$cl$Q_{V}:V\in $\ss $\},$
then for each $V\in \text{\ss },$ $x^{\ast }\in $cl$\left\{ x^{\ast }\in
K:x^{\ast }\in (T(x^{\ast })+V)\cap K\right\} $. Hence $(x^{\ast },x^{\ast
})\in $clGr($(T(x)+V)\cap K)$ for each $V\in $\ss . By Lemma 2\emph{\ }we
have that\emph{\ } $x^{\ast }\in \overline{T}(x^{\ast }),$ i.e. $x^{\ast }$
is a fixed point for $\overline{T}$ $\Box $\medskip

The weakly convex correspondences are defined in [7].

\begin{definition}
\textit{\ }[7].\textit{\ }Let $X$ and $Y$ be nonempty convex subsets of a
topological vector space $E.$\textit{\ }The correspondence\textit{\ }$%
T:X\longrightarrow 2^{Y}$ is said to have \emph{weakly convex graph} (in
short it is a WCG correspondence), if for each finite set $%
\{x_{1},x_{2},...,x_{n}\}\subset X$, there exists $y_{i}\in T(x_{i})$, $%
(i=1,2,...,n)$, such that
\end{definition}

\begin{center}
(1) $\ $co$(\{(x_{1},y_{1}),(x_{2},y_{2}),...,(x_{n},y_{n})\})\subset $Gr$%
(T) $
\end{center}

\QTP{Body Math}
The relation (1) is equivalent to

\begin{center}
(2) $\ \ \ \ \overset{n}{\underset{i=1}{\tsum }}\lambda _{i}y_{i}\in T(%
\overset{n}{\underset{i=1}{\tsum }}\lambda _{i}x_{i})$ \ \ \ \ \ \ \ $%
(\forall (\lambda _{1},\lambda _{2},...,\lambda _{n})\in \Delta _{n-1}).$
\end{center}

It is clear that if either Gr$(T)$ is convex, or $\ \tbigcap \{T(x):x\in
X\}\neq \emptyset $, then $T$ has a weakly convex graph.\medskip

\textit{Remark} 4. Let $T:X\rightarrow 2^{Y}$ \ be a WCG correspondence and $%
X_{0}$ be a non-empty convex subset of $X$. Then, the restriction of $T$ on $%
X_{0}$, $T_{\mid X_{0}}:X_{0}\rightarrow 2^{Y}$ is a WCG correspondence,
too.\medskip

Now we introduce the following definition.

\begin{definition}
Let $E$, $F$ be topological vector spaces, $X$ and $Y$ be nonempty convex
subsets of $E,$ respectively $F$ and $T:X\rightarrow 2^{Y}$ be a
correspondence. $T$ is said to have \textit{a *-weakly convex graph} if for
each neighborhood $V$ of the origin in $F,$ the correspondence $%
T_{V}:X\rightarrow 2^{Y},$ defined by $T_{V}(x)=(T(x)+V)\cap Y$ for each $%
x\in X$ has an weakly convex graph.
\end{definition}

The next theorem (its proof follows the same lines as that of Theorem 5) is
a fixed point result for a correspondence with *-weakly convex graph.

\begin{theorem}
Let $K$ be a simplex in a topological vector space $F.$ Let $T:K\rightarrow
2^{K}$ be a correspondence with *-weakly convex graph. Then, there exists a
point $x^{\ast }\in K$ such that $x^{\ast }\in \overline{T}(x^{\ast }).$
\end{theorem}

We get the following corollary$.$

\begin{corollary}
Let $K$ be a simplex in a topological vector space $F.$ Let \textit{\ }$%
S,T:K\rightarrow 2^{K}$\textit{\ be two correspondences with the following
conditions:}
\end{corollary}

(i) \textit{for each }$x\in K,$\textit{\ }$\overline{S}(x)\subset T(x)$%
\textit{\ and }$S(x)\neq \emptyset ,$

(ii)\textit{\ }$S$\textit{\ has *-weakly convex graph.}

\textit{Then, there exists a point }$x^{\ast }\in K$\textit{\ such that }$%
x^{\ast }\in T(x^{\ast }).\medskip $

Now, we introduce the concept of *-weakly naturally quasiconvex
correspondence.

\begin{definition}
Let $E$, $F$ be topological vector spaces, $X$ and $Y$ be nonempty convex
subsets of $E,$ respectively $F$ and $T:X\rightarrow 2^{Y}$ be a
correspondence. $T$ is said to be\textit{\ *-weakly naturally quasiconvex}
if for each neighborhood $V$ of the origin in $F,$ the corespondence $%
T_{V}:X\rightarrow 2^{Y},$ defined by $T_{V}(x)=(T(x)+V)\cap Y$ for each $%
x\in X$ is weakly naturally quasiconvex.
\end{definition}

Theorem 7 is a fixed point theorem for \textit{*-}weakly naturally
quasiconvex\textit{\ }correspondences.

\begin{theorem}
Let $K$ be a non-empty simplex in a topological vector space $F.$ Let $%
T:K\rightarrow 2^{K}$ be a \textit{*-weakly naturally quasiconvex }%
correspondence. Then there exists a point $x^{\ast }\in K$ such that $%
x^{\ast }\in \overline{T}(x^{\ast }).$
\end{theorem}

\textit{Proof.} Let \ss\ denote the family of all neighborhoods of zero in $%
F $ and let $V\in $\ss $.$ The corespondence $T_{V}:K\rightarrow 2^{K},$
defined by $T_{V}(x)=(T(x)+V)\cap K$ for each $x\in K$ is *-weakly naturally
quasiconvex. Then there exists a continuous selection $f_{V}:K\rightarrow K$
such that $f_{V}(x)\in T_{V}(x).$ The proof follows the same line as in
Theorem 5.$\Box $\medskip

Now we introduce the following definition.

\begin{definition}
Let $B\subset X\times Y$ be a biconvex set, $Z$ a nonempty convex subset of
a topological vector space $F$ and $T:B\rightarrow 2^{Z}$ a correspondence$.$
$T$ is called \textit{weakly biconvex} if for each finite set $%
\{(x_{1},y_{1}),(x_{2},y_{2}),...,(x_{n},y_{n})\}\subset B$, there exists $%
z_{i}\in T(x_{i},y_{i})$, $(i=1,2,...,n)$ such that for every biconvex
combination $(x,y)=\tsum\limits_{i=1}^{n}\lambda _{i}(x_{i},y_{i})\in B$
(with $\tsum\limits_{i=1}^{n}\lambda _{i}=1,$ $\lambda _{i}\geq 0$ $%
i=1,2,...,n$)$,$ there exists $y^{\prime }\in T(\overset{n}{\underset{i=1}{%
\tsum }}\lambda _{i}(x_{i},y_{i}))$ and $y^{\prime }=\overset{n}{\underset{%
i=1}{\tsum }}\lambda _{i}z_{i}.\medskip $
\end{definition}

We state the following selection theorem for weakly biconvex
correspondences.\medskip

\begin{theorem}
(selection theorem). \textit{Let }$Y$\textit{\ be a non-empty convex subset
of a topological vector space }$F$\textit{\ and }$K\subset E_{1}\times
E_{2}, $ where $E_{1},E_{2}$ are topological vector spaces$.$ Suppose that $%
K $ is the biconvex hull of $\{(a_{1},b_{1}),(a_{2},b_{2}),...,(a_{n},b_{n})%
\}\subset E_{1}\times E_{2}$\textit{.} \textit{\ Let }$T:K\rightarrow 2^{Y}$%
\textit{\ be a weakly biconvex correspondence. Then, }$T$\textit{\ has a
continuous selection on }$K$\textit{.\medskip }
\end{theorem}

\textit{Proof.} Since $T$ is weakly biconvex, there exists $c_{i}\in
T(a_{i},b_{i})$, $(i=1,2,...,n),$ such that for every $(\lambda _{1},\lambda
_{2},...,\lambda _{n})\in \Delta _{n-1}$, there exists $z\in T(\overset{n}{%
\underset{i=1}{\tsum }}\lambda _{i}(a_{i},b_{i}))$ with $z=\overset{n}{%
\underset{i=1}{\tsum }}\lambda _{i}z_{i}.$

Since $K$ is biconvex hull of $(a_{1},b_{1}),...,(a_{n},b_{n}),$ there exist
unique continuous functions $\lambda _{i}:K\rightarrow \mathbb{R},$ $%
i=1,2,...,n$ such that for each $(x,y)\in K,$ we have $(\lambda
_{1}(x,y),\lambda _{2}(x,y),...,\lambda _{n}(x,y))\in \Delta _{n-1}$ and $%
(x,y)=\overset{n}{\underset{i=1}{\tsum }}\lambda _{i}(x,y)(a_{i},b_{i}).$

Define $f:K\rightarrow 2^{Y}$ by

$f(a_{i},b_{i})=c_{i}$ $(i=1,...,n)$ and

$f(\overset{n}{\underset{i=1}{\tsum }}\lambda _{i}(a_{i},b_{i}))=\overset{n}{%
\underset{i=1}{\tsum }}\lambda _{i}c_{i}\in T(x,y).$

We show that $f$ is continuous.

Let $(x_{m},y_{m})_{m\in N}$ be a sequence which converges to $x_{0}\in K,$
where $(x_{m},y_{m})=\overset{n}{\underset{i=1}{\tsum }}\lambda
_{i}(x_{m},y_{m})(a_{i},b_{i})$ implies $a_{1}=a_{2}=...=a_{n}=a$ or $%
b_{1}=b_{2}=...=b_{n}=b$ and $(x_{0},y_{0})=$ $\overset{n}{\underset{i=1}{%
\tsum }}\lambda _{i}(x_{0})(a_{i},b_{i})$ with $a_{1}=a_{2}=...=a_{n}=a$ or $%
b_{1}=b_{2}=...=b_{n}=b.$ By the continuity of $\lambda _{i},$ it follows
that for each $i=1,2,...,n$, $\lambda _{i}(x_{m},y_{m})\rightarrow \lambda
_{i}(x_{0},y_{0})$ as $m\rightarrow \infty .$ Hence $f(x_{m},y_{m})%
\rightarrow f(x_{0},y_{0})$ as $m\rightarrow \infty ,$ i.e. $f$ is
continuous.

We proved that $T$ has a continuous selection on $K.$\medskip

In order to prove the existence theorems of equilibria for a generalized
quasi-game, we need the following version of Kim's quasi fixed-point theorem:

\smallskip

\begin{theorem}
\textit{Let }$I$\textit{\ and }$J$\textit{\ be any (possible uncountable)
index sets. For each }$i\in I$\textit{\ and }$j\in J$\textit{, let }$X_{i}$%
\textit{\ and }$Y_{j}$\textit{\ be non-empty compact convex subsets of
Hausdorff locally convex spaces }$E_{i}$\textit{\ and respectively }$F_{j}$%
\textit{. }
\end{theorem}

\textit{Let }$X:=\prod X_{i}$\textit{, }$Y:=\underset{i\in I}{\prod }Y_{j}$%
\textit{\ and }$Z:=X\times Y$\textit{.}

\textit{For each }$i\in I$\textit{\ let }$\Phi _{i}:Z\rightarrow 2^{X_{i}}$%
\textit{\ be a correspondence such that the set }$W_{i}=\left\{ (x,y)\in Z%
\text{ }\mid \Phi _{i}(x,y)\neq \emptyset \right\} $\textit{\ is open and }$%
\Phi _{i}$\textit{\ has a continuous selection f}$_{i}$ \textit{on }$W_{i}$%
\textit{.}

\textit{For each j}$\in J$\textit{\ let }$\Psi _{j}:Z\rightarrow 2^{Y_{j}}$%
\textit{\ be an upper semicontinuous correspondence with non-empty closed
convex values.}

\textit{Then there exists a point }$(x^{\ast },y^{\ast })\in Z$\textit{\
such that for each }$i\in I$\textit{, either }$\Phi _{i}(x^{\ast },y^{\ast
})=\emptyset $\textit{\ or }$\overset{\_}{x}_{i}\in \Phi _{i}(x^{\ast
},y^{\ast })$\textit{, and for each }$j\in J$\textit{, }$y_{j}^{\ast }\in
\Psi _{j}(x^{\ast },y^{\ast })$\textit{.\medskip }

\textit{\ }\smallskip \textit{Proof.} We first endow $\underset{i\in I}{%
\prod }E_{i}$ and $\underset{j\in J}{\prod }F_{j}$ with the product
topologies; and then $\underset{i\in I}{\prod }E_{i}\times $ $\underset{j\in
J}{\prod }F_{j}$ is also a locally convex Hausdorff topological vector space.

For each $i\in I$, we define a correspondence $\Phi _{i}^{^{\prime
}}:Z\rightarrow 2^{X_{i}}$ by

$\Phi _{i}^{^{\prime }}(x,y):=\left\{ 
\begin{array}{c}
\{f_{i}(x,y)\}\text{, if }(x,y)\in W_{i}\text{, } \\ 
X_{i}\text{, \qquad if }(x,y)\notin W_{i}\text{;}%
\end{array}%
\right. $

Then for each $(x,y)\in Z$, $\Phi _{i}^{^{\prime }}(x,y)$ is a non-empty
closed convex subset of $X_{i}$. Also, $\Phi _{i}^{^{\prime }}$ is an upper
semicontinuous correspondence on $Z$. In fact, for each proper open subset $%
V $ of $X_{i}$, we have

$U:=\{(x,y)\in Z$ $\mid $ $\Phi _{i}^{^{\prime }}(x,y)\subset V\}$

\ \ \ =$\{(x,y)\in W_{i}$ $\mid $ $\Phi _{i}^{^{\prime }}(x,y)\subset V\cup
(x,y)\in Z\setminus W_{i}$ $\mid $ $\Phi _{i}^{^{\prime }}(x,y)\subset V\}$

...=$\left\{ (x,y)\in W_{i}\text{ }\mid \text{ }f_{i}(x,y)\in V\right\} \cup
\left\{ (x,y)\in Z\setminus W_{i}\text{ }\mid \text{ }X_{i}\subset V\right\} 
$

\ \ \ =$\left\{ (x,y)\in W_{i}\text{ }\mid \text{ }f_{i}(x,y)\in V\right\}
=f_{i}^{-1}(V)\cap W_{i}$.

Since $W_{i}$ is open and $f_{i}$ is a continuous map on $W_{i},U$ is open,
and hence $\Phi _{i}^{^{\prime }}$ is upper semicontinuous on $Z$.

Finally, we define a correspondence $\Phi :Z\rightarrow 2^{Z}$ by

$\Phi (x,y):=\underset{i\in I}{\prod }\Phi _{i}^{^{\prime }}(x,y)\times 
\underset{j\in J}{\prod }\Psi _{j}(x,y)$ for each $(x,y)\in Z$.

Then, by Lemma 3 in [8], $\Phi $ is an upper semicontinuous correspondence
such that each $\Phi (x,y)$ is non-empty closed convex. Therefore, by the
Fan-Glicksebrg fixed point theorem [9] there exists a fixed point $(x^{\ast
} $,$y^{\ast })\in Z$ such that $(x^{\ast },y^{\ast })$ $\in \Phi (x,y)$,
i.e., for each $i\in I$, $x_{i}^{\ast }\in \Phi _{i}^{^{\prime }}(x,y)$, and
for each $j\in J$, $y_{j}^{\ast }\in \Psi _{j}(x,y)$. If $(x^{\ast },y^{\ast
})\in W_{i}$ for some $i\in I$, then $x_{i}^{\ast }=f_{i}(x^{\ast },y^{\ast
})\in \Phi _{i}(x^{\ast },y^{\ast })$; and if $(x^{\ast },y^{\ast })\notin
W_{i}$ for some $i\in I$, then $\Phi _{i}(x^{\ast },y^{\ast })=\emptyset $.
Therefore, we have that for each $i\in I$, either $\Phi _{i}(x^{\ast
},y^{\ast })=\emptyset $ or $x_{i}^{\ast }\in \Phi _{i}(x^{\ast },y^{\ast })$%
. Also, for each $j\in J$, we already have $y_{j}^{\ast }\in \Psi
_{j}(x^{\ast },y^{\ast })$. This completes the proofs. $\Box \medskip $

We have the following corollary.

\begin{corollary}
\textit{Let }$I$\textit{\ and }$J$\textit{\ be any (possible uncountable)
index sets. For each }$i\in I$\textit{\ and }$j\in J$\textit{, let }$X_{i}$%
\textit{\ and }$Y_{j}$\textit{\ be non-empty compact convex subsets of
Hausdorff locally convex spaces }$E_{i}$\textit{\ and respectivelly }$F_{j}$%
\textit{. }
\end{corollary}

\textit{Let }$X:=\prod X_{i}$\textit{, }$Y:=\underset{i\in I}{\prod }Y_{j}$%
\textit{\ and }$Z:=X\times Y$\textit{.}

\textit{For each }$i\in I$\textit{\ let }$S_{i}:Z\rightarrow 2^{X_{i}}$%
\textit{\ be a correspondence such that the set }$W_{i}=\left\{ (x,y)\in Z%
\text{ }\mid S_{i}(x,y)\neq \emptyset \right\} $\textit{\ is the interior of
the biconvex hull of }$\{(a_{1},b_{1}),(a_{2},b_{2}),...,(a_{n},b_{n})\}%
\subset Z$ \textit{and }$S_{i}$\textit{\ is weakly biconvex} \textit{on }$%
W_{i}$\textit{.}

\textit{For each }$\mathit{j}\in J$\textit{\ let }$T_{j}:Z\rightarrow
2^{Y_{j}}$\textit{\ be an upper semicontinuous correspondence with non-empty
closed convex values.}

\textit{Then there exists a point }$(x^{\ast },y^{\ast })\in Z$\textit{\
such that for each }$i\in I$\textit{, either }$S_{i}(x^{\ast },y^{\ast
})=\emptyset $\textit{\ or }$x_{i}^{\ast }\in S_{i}(x^{\ast },y^{\ast })$%
\textit{, and for each }$j\in J$\textit{, }$y_{j}^{\ast }\in T_{j}(x^{\ast
},y^{\ast })$\textit{.\medskip }

\section{Applications in the equilibrium theory\textbf{\protect\medskip }}

In this paper, we study the following model of a generalized
quasi-game.\medskip

\begin{definition}
\textit{\ }Let $I$ be a nonempty set (the set of agents). For each $i\in I$,
let $X_{i}$ be a non-empty topological vector space representing the set of
actions and define $X:=\underset{i\in I}{\prod }X_{i}$; let $A_{i}$, $%
B_{i}:X\times X\rightarrow 2^{X_{i}}$ be the constraint correspondences and $%
P_{i}:X\times X\rightarrow 2^{X_{i}}$ the preference correspondence. A 
\textit{generalized quasi-game} $\Gamma =(X_{i},A_{i},B_{i},P_{i})_{i\in I}$
is defined as a family of ordered quadruples $(X_{i},A_{i},B_{i},P_{i})$.
\end{definition}

In particular, when $I$=$\left\{ 1\text{, }2\text{...}n\right\} $, $\Gamma $
is called n-person quasi-game.\medskip

\begin{definition}
\textit{\ }An \textit{equilibrium} for $\Gamma $ is defined as a point $%
(x^{\ast },y^{\ast })\in X\times X$ such that for each $i\in I$, $%
y_{i}^{\ast }\in $cl$B_{i}(x^{\ast },y^{\ast })$ and $A_{i}(x^{\ast
},y^{\ast })\cap P_{i}(x^{\ast },y^{\ast })=\emptyset $.\medskip
\end{definition}

If $A_{i}(x,y)=B_{i}(x,y)$ for each $(x,y)\in X\times X$ and $i\in I$, this
model coincides with that one introduced by W. K. Kim [13].

If, in addition, for each $i\in I$, $\ A_{i},P_{i}$ are constant with
respect to the first argument, this model coincides with the classical one
of the abstract economy and the definition of equilibrium is that given in
[4].

In this work, Kim established an existence result for a generalized
quasi-game with a possibly uncountable set of agents, in a locally convex
Hausdorff topological vector space.

\smallskip Here is his result:\medskip

\begin{theorem}[11]
\textit{Let }$\Gamma =(X_{i},A_{i},B_{i},P_{i})_{i\in I}$\textit{\ be a
generalized quasi-game, where }$I$\textit{\ is a (possibly uncountable) set
of agents such that for each }$i\in I$\textit{\ :}
\end{theorem}

\textit{(1) }$X_{i}$\textit{\ is a non-empty compact convex subset of a
Hausdorff \ locally convex space }$E_{i}$\textit{\ and denote }$X:=\underset{%
i\in I}{\prod }X_{i}$\textit{\ and }$Z:=X\times X$\textit{;}

\textit{(2) The correspondence }$A_{i}:X\times X\rightarrow 2^{X_{i}}$%
\textit{\ is upper semicontinuous such that }$A_{i}(x,y)$ \textit{is a
non-empty convex subset of }$X_{i}$\textit{\ for each }$(x,y)\in Z$\textit{;}

\textit{(3) }$A_{i}^{-1}\left( x_{i}\right) $\textit{\ is (possibly empty)
open for each }$x_{i}\in X_{i}$\textit{;}

\textit{(4) the correspondence }$P_{i}:Z\rightarrow 2^{X_{i}}$\textit{\ is
such that }$(A_{i}\cap P_{i})^{-1}\left( x_{i}\right) $\textit{\ is
(possibly empty) open for each }$x_{i}\in X_{i}$\textit{;}

\textit{(5) the set }$W_{i}$\textit{\ }$:$\textit{\ }$=\left\{ (x,y)\in Z%
\text{ }\mid \text{ }\left( A_{i}\cap P_{i}\right) (x,y)\neq \emptyset
\right\} $\textit{\ is perfectly normal;}

\textit{(6) for each }$(x,y)\in W_{i}$\textit{, }$x_{i}\notin coP_{i}(x,y)$%
\textit{.}

\textit{Then there exists an equilibrium point }$(x^{\ast },y^{\ast })\in
X\times X$\textit{\ for }$\Gamma $\textit{, }$i.e$\textit{., for each }$i\in
I$\textit{, }$y_{i}^{\ast }\in $cl$A_{i}(x^{\ast },y^{\ast })$\textit{\ and }%
$A_{i}(x^{\ast },y^{\ast })\cap P_{i}(x^{\ast },y^{\ast })=\emptyset $%
\textit{.\medskip }

As application of the selection theorems from section 3, we state a theorem
on the existence of the equilibrium for a generalized quasi-game.\medskip

\begin{theorem}
\ \textit{Let }$\Gamma =(X_{i},A_{i},B_{i},P_{i})_{i\in I}$\textit{\ \ be a
generalized quasi-game where }$I$\textit{\ is a (possibly uncountable) set
of agents such that for each }$i\in I:$
\end{theorem}

\textit{(1) }$X_{i}$\textit{\ is a non-empty compact convex set in a
Hausdorff locally convex space }$E_{i}$\textit{\ and denote }$X:=\underset{%
i\in I}{\prod }X_{i}$\textit{\ and }$Z:=X\times X$\textit{;}

\textit{(2) The correspondence }$B_{i}:Z\rightarrow 2^{X_{i}}$\textit{\ is
non-empty, convex valued such that for each }($x,y)\in Z$,\textit{\ }$%
A_{i}(x,y)\subset B_{i}(x,y)\ $\textit{and }cl$B_{i}$\textit{\ is upper
semicontinuous;}

\textit{(3) the correspondence }$A_{i}\cap P_{i}:W_{i}\rightarrow 2^{X_{i}}$%
\textit{\ is weakly naturally quasiconvex;}

\textit{(4) the set }$W_{i}:$\textit{\ }$=\left\{ (x,y)\in Z\text{ / }\left(
A_{i}\cap P_{i}\right) (x,y)\neq \emptyset \right\} $\textit{\ is open and }%
cl$W_{i}$ \textit{is a }$(n-1)$\textit{\ dimensional simplex in }$Z$;

\textit{(5) for each }$(x,y)\in W_{i},$\textit{\ }$x_{i}\notin P_{i}(x,y)$%
\textit{.}

\textit{Then there exists an equilibrium point }$(x^{\ast },y^{\ast })\in Z$ 
\textit{\ for }$\Gamma $\textit{,}$\ i.e.$\textit{, for each }$i\in I$%
\textit{, }$y_{i}^{\ast }\in $cl$B_{i}(x^{\ast },y^{\ast })$\textit{\ and }$%
A_{i}(x^{\ast },y^{\ast })\cap P_{i}(x^{\ast },y^{\ast })=\emptyset $\textit{%
.\bigskip }

\textit{Proof.} For each $i\in I$, we define $\Phi _{i}:Z\rightarrow
2^{X_{i}}$ by

$\Phi _{i}(x,y)=\left\{ 
\begin{array}{c}
(A_{i}\cap P_{i})(x,y)\text{, if }(x,y)\in W_{i}\text{, } \\ 
\emptyset \text{, \qquad \qquad \qquad\ \ if }(x,y)\notin W_{i}\text{;}%
\end{array}%
\right. $

By applying Theorem 3 to the restrictions $A_{i}\cap P_{i}$ on $W_{i}$, we
can obtain that there exists a continuous selection $f_{i}:W_{i}\rightarrow
X_{i}$ such that $f_{i}(x,y)\in (A_{i}\cap P_{i})(x,y)$ for each $(x,y)\in
W_{i}$.

For each $j\in I$, we define $\Psi _{j}:Z\rightarrow 2^{X_{i}}$, by $\Psi
_{j}(x,y)=$cl$B_{j}(x,y)$ for each $(x,y)\in Z$.

Then $\Psi _{j}$ is an upper semicontinuous correspondence and $\Psi
_{j}(x,y)$ is a non-empty, convex, closed subset of $X_{j}$ for each $%
(x,y)\in Z$.

By Theorem 9, it follows that there exists $(x^{\ast },y^{\ast })\in Z$ such
that for each $i\in I$, either $\Phi _{i}(x^{\ast },y^{\ast })=\emptyset $
or $x_{i}^{\ast }\in \Phi _{i}(x^{\ast },y^{\ast })$ and for each $j\in J$, $%
y_{j}^{\ast }\in \Psi _{j}(x^{\ast },y^{\ast })$.

If $x_{i}^{\ast }\in \Phi _{i}(x^{\ast },y^{\ast })$ for some $i\in I$, then 
$x_{i}^{\ast }\in \Phi _{i}(x^{\ast },y^{\ast })=(A_{i}\cap P_{i})(x^{\ast
},y^{\ast })\subset P_{i}(x^{\ast },y^{\ast })$ which contradicts the
assumption (5).

Therefore, for each $i\in I$, $\Phi _{i}(x,y)=\emptyset $ and then $(x^{\ast
},y^{\ast })\notin W_{i}$. Hence, $(A_{i}\cap P_{i})(x^{\ast },y^{\ast
})=\emptyset $ and for each $i\in I$, $y^{\ast }\in \Psi _{i}(x^{\ast
},y^{\ast })=$cl$B_{i}(x^{\ast },y^{\ast })$. $\Box $

By using a similar type of proof and Theorem 8, we obtain Theorem 12.

\begin{theorem}
\textit{Let }$\Gamma =(X_{i},A_{i},B_{i},P_{i})_{i\in I}$\textit{\ \ be a
generalized quasi-game where }$I$\textit{\ is a (possibly uncountable) set
of agents such that for each }$i\in I:$
\end{theorem}

\textit{(1) }$X_{i}$\textit{\ is a non-empty compact convex set in a
Hausdorff locally convex space }$E_{i}$\textit{\ and denote }$X:=\underset{%
i\in I}{\prod }X_{i}$\textit{\ and }$Z:=X\times X$\textit{;}

\textit{(2) The correspondence }$B_{i}:Z\rightarrow 2^{X_{i}}$\textit{\ is
non-empty, convex valued such that for each }($x,y)\in Z$,\textit{\ }$%
A_{i}(x,y)\subset B_{i}(x,y)\ $\textit{and }cl$B_{i}$\textit{\ is upper
semicontinuous;}

\textit{(3) }$A_{i}\cap P_{i}$\textit{\ is a weakly biconvex correspondence
on }$W_{i}$\textit{;}

\textit{(4) the set }$W_{i}:$\textit{\ }$=\left\{ (x,y)\in Z\text{ / }\left(
A_{i}\cap P_{i}\right) (x,y)\neq \emptyset \right\} $\textit{\ is the
interior of the biconvex hull of }$%
\{(a_{1},b_{1}),(a_{2},b_{2}),...,(a_{n},b_{n})\}\subset Z$\textit{;}

\textit{(5) for each }$(x,y)\in W_{i},$\textit{\ }$x_{i}\notin $co$%
P_{i}(x,y) $\textit{.}

\textit{Then there exists an equilibrium point }$(x^{\ast },y^{\ast })\in Z$ 
\textit{\ for }$\Gamma $\textit{,}$\ i.e.$\textit{, for each }$i\in I$%
\textit{, }$y_{i}^{\ast }\in $cl$B_{i}(x^{\ast },y^{\ast })$\textit{\ and }$%
A_{i}(x^{\ast },y^{\ast })\cap P_{i}(x^{\ast },y^{\ast })=\emptyset $\textit{%
.\medskip }


\begin{thebibliography}{99}
\bibitem{1} R. P. Agarwal, O'Regan, \textit{A Note on Equilibria for
abstract economies,} Mathematical and Computer Modelling 34 (2001), 331-343.

\bibitem{2} R. Aumann, S. Hart, \textit{Bi-convexity and bi-martingales,}
Isr J Math \textbf{54, }\textit{2} (1986), 159-180.

\bibitem{3} K. J. Arrow, G. Debreu, \textit{Existence of an Equilibrium for
a Competitive} \textit{Economy}. Econometrica \textbf{22} (1954), 265-290.

\bibitem{4} A. Borglin and H. Keiding, \textit{Existence of equilibrium
action and of equilibrium:A note on the "new" existence theorem.} J. Math.
Econom. \textbf{3} (1976), 313-316.

\bibitem{5} F. E. Browder, \textit{A new generation of the Schauder fixed
point theorems.} Math. Ann. \textbf{174} (1967), 285-290.

\bibitem{6} G. Debreu, \textit{A social equilibrium existence theorem.}
Proc. Nat. Acad. Sci. \textbf{38} (1952), 886-893.

\bibitem{7} X. Ding, He Yiran, \textit{Best Approximation Theorem for
Set-valued Mappings without Convex Values and} \textit{Continuity.} Appl
Math. and Mech. English Edition, \textbf{19}, \textit{9} (1998), 831-836.

\bibitem{8} K. Fan, \textit{Fixed-point and minimax theorems in locally
convex topological linear spaces. }Proc. Nat. Acad. Sci. U.S.A. \textbf{38 }%
(1952), 121-126.

\bibitem{9} K. Fan, \textit{A generalization of Tychonoff's fixed point
theorem}, Math. Ann. 142 (1961), 305-310.

\bibitem{10} J. Gorski, F. Pfeuffer, K. Klamroth, \textit{Biconvex sets and
optimization with biconvex functions: a survey and extension, }Math. Meth.
Oper. Res. 66 (2007), 373-407.

\bibitem{11} C.D. Horvath, \textit{Point fixes et coincidences pour les
applications multivoques sans convexite}, C. R. Acad. Sci. Paris 296 (1983),
403-406.

\bibitem{12} W. K. Kim, K. K. Tan, \textit{New existence theorems of
equilibria and applications.} Nonlinear Anal. \textbf{47} (2001), 531-542.

\bibitem{13} W. K. Kim, \textit{On a quasi fixed-point theorem. }Bull.
Korean Math. Soc. \textbf{40} (2003), \textit{2}, 301-307.

\bibitem{14} I. Lefebvre, \textit{A quasi fixed-point theorem for a product
of u.s.c. or l.s.c. correspondences with an economic application,}
Set-Valued Anal. \textbf{9} (2001), 273-288.

\bibitem{15} J.F. Nash, \textit{Non-cooperative games.} Ann. Math. \textbf{54%
} (1951), 286-295.

\bibitem{16} W. Shafer, H. Sonnenschein, \textit{Equilibrium in abstract
economies without ordered} \textit{preferences. }J. Math. Econom. \textbf{2}
(1975), 345-348.

\bibitem{17} G. Tian, \textit{Fixed points theorems for mappings with
non-compact and non-convex domains,} J. Math. An. and Appl. \textbf{158}
(1991), 161-167.

\bibitem{18} K. Vind, \textit{Equilibrium with coordination,} J. Math.
Econom. \textbf{12} (1983), 275-285.

\bibitem{19} N. C. Yannelis and N. D. Prabhakar, \textit{Existence of
maximal elements and equilibrium} \textit{in linear topological spaces.} J.
Math. Econom. \textbf{12} (1983), 233-245.

\bibitem{20} X. Z. Yuan, \textit{The Study of Minimax inequalities and
Applications to Economies and} \textit{Variational inequalities. }Memoirs of
the American Society \textbf{132}, \textit{625, }(1988).
\end{thebibliography}
\end{document}